\documentclass{ifacconf}

\usepackage{graphicx}      
\usepackage{natbib}        
\usepackage{amsmath}
\usepackage{amssymb}
\usepackage{amsfonts}
\usepackage{mathrsfs}
\usepackage{bm}
\usepackage{physics}
\usepackage{mathtools}
\usepackage{algorithm}
\usepackage{algpseudocode}
\usepackage{booktabs}
\usepackage{siunitx}
\usepackage{subcaption}
\usepackage{booktabs}
\begin{document}
\begin{frontmatter}

\title{Trajectory Optimization by \\Pseudospectral Successive Convexification on Riemannian Manifolds} 

\thanks[footnoteinfo]{\copyright~2026 the authors.
This work has been accepted to IFAC for publication under a
Creative Commons Licence CC-BY-NC-ND.}
\author[First]{Tatsuya Narumi} 
\author[Second]{Shin{-}ichiro Sakai}

\address[First]{Department of Advanced Energy, 
   The University of Tokyo, Chiba, Japan (e-mail: narumi.tatsuya@ac.jaxa.jp).}
\address[Second]{JAXA Institute of Space and Astronautical Science, 
   Kanagawa, Japan (e-mail: sakai@isas.jaxa.jp).}

\begin{abstract}
This paper proposes an intrinsic pseudospectral convexification framework 
for optimal control problems with manifold constraints.
While pseudospectral successive  convexification combines spectral collocation 
with successive convexification, classical pseudospectral methods are 
not geometry-consistent on manifolds. 
This is because interpolation and differentiation are performed in Euclidean coordinates.
We introduce a geometry-consistent transcription that enables pseudospectral collocation 
without imposing manifold constraints extrinsically. 
The resulting method solves nonconvex manifold-constrained problems 
through a sequence of convex subproblems. 
A six-degree-of-freedom landing guidance example with unit quaternions and unit direction vectors demonstrates the practicality of the approach. 
The proposed method preserves manifold feasibility to machine precision  and achieves significant computational speedups. 

\end{abstract}


\end{frontmatter}

\section{Introduction}

A central challenge in trajectory optimization is to faithfully represent the complex, 
nonconvex nature of physical systems within a computationally tractable framework.
To address this challenge, Successive Convexification (SCvx) by \cite{mao_scvx}, has emerged as a predominant direct optimization approach.
This framework iteratively solves convex approximations of the original nonconvex problem, constructed by linearizing dynamics and constraints around a reference trajectory.
Owing to its tractability and favorable convergence properties, SCvx has become a standard tool in optimal control.
Recently, its scope has been expanded to handle even logical complexity; for instance, \cite{szmuk_stc}, demonstrated that if/then constraints 
can be incorporated within the SCvx framework.

Despite this progress in SCvx, 
the treatment of \textit{geometric} complexity---
specifically, nonlinear manifold constraints such as the quaternion unit-norm constraint $\norm{q} = 1$---
remains neglected (e.g. \cite{szmuk_Free6DoF}) or is handled extrinsically.
Common extrinsic approaches are Baumgarte stabilization by \cite{wojciech_stabilization}
and convex-concave decomposition (CCD) by \cite{lu_CCD} and \cite{sagliano_accD}.
They embed the manifold in Euclidean space as a nonlinear equality $h(x) = 0$.
However, these extrinsic strategies do not strictly respect the manifold structure 
and incur computational overhead, 
typically requiring auxiliary constraints or variables 
for a single scalar constraint.
Crucially, Baumgarte stabilization may cause energy dispersion and 
instabilities, while CCD's applicability is limited to 
manifolds defined by $h(x) = 0$, where $h$ is convex. 

Intrinsic methods, in contrast, preserve the geometry explicitly 
through operations on tangent spaces and retractions. 
In particular, \cite{kraisler_iscvx} proposed an intrinsic SCvx framework 
that avoids reliance on a specific Euclidean embedding and supports a wide range 
of Riemannian manifolds via retraction-based local parameterizations. 
However, this intrinsic formulation is discrete-time (e.g., $x_{k+1}=f(x_k,u_k)$) and 
therefore does not directly interface with pseudospectral (PS) transcriptions by \cite{sagliano_SPsCvx}, 
where global interpolation and differentiation enable \emph{spectral convergence} 
(i.e., high accuracy with few nodes). 
This gap prevents combining intrinsic SCvx with the numerical efficiency of spectral collocation.

Motivated by tangent-space collocation formulations 
for direct methods \cite{Bordalba_direct_collocation} and 
Lie group optimization \cite{Saccon_LieGroupProjection}, 
we develop an intrinsic  pseudospectral SCvx framework for manifold-constrained optimal control.
The main contributions of this paper are 
\begin{itemize}
\item We propose an \emph{intrinsic pseudospectral (iPS) method} in which state and input trajectories are parameterized by retractions, preserving manifold feasibility by construction without imposing extrinsic constraints.
\item We propose a \emph{manifold-consistent collocation scheme}. 
By transporting nodal perturbations to the same tangent space, this scheme enables the standard application of PS methods and their integration with SCvx.
\item We demonstrate the resulting algorithm on a six-degree-of-freedom powered-landing problem with unit quaternions and unit thrust-direction vectors.
The proposed method shows comparable optimality to an extrinsic ACCD baseline while eliminating unit-norm drift, and achieves significant computational speedups due to reduced problem dimensionality.
\end{itemize}

The remainder of this paper is organized as follows. Section~2 
reviews preliminaries on Riemannian manifolds. Sections~3 and~4 formulate the problem and discuss the 
limitations of standard PS methods. Section~5 
presents the proposed iPS framework, which is validated
 via a six-degree-of-freedom landing simulation in Section~6. Section~7 concludes the paper.

\section{Preliminaries on Riemannian manifolds}
In this section, we summarize the basic geometric operators.
See \cite{Absil_optManifold} and \cite{boumal2023intromanifolds}, for details.

Let $\mathcal M$ denote a smooth manifold, which locally resembles some Euclidean space $\mathbb R^n$.  
and $T_x\mathcal M$ its tangent space at $x\in\mathcal M$, which is a vector space consisting of tangent vectors at $x$.
The tangent bundle $T\mathcal M$ is the collection of all tangent spaces.

\textbf{Directional derivative:}
For a smooth map $f:\mathcal M\to\mathcal N$, the directional derivative along
$\eta\in T_x\mathcal M$ is defined by
\begin{align}
  Df(x)[\eta]\coloneqq \frac{d}{dt} f(\gamma(t))\Big|_{t=0},
  \quad \gamma(0)=x,\ \dot\gamma(0)=\eta.
\end{align}

\textbf{Retraction:}
A retraction is a smooth map $R_x:T_x\mathcal M\to\mathcal M$ satisfying
\begin{align}\label{eq:retraction_def}
  R_x(0_x)=x,\qquad DR_x(0_x)=\mathrm{id}_{T_x\mathcal M}.
\end{align}
By this retraction, a neighborhood of $x\in\mathcal M$ is parameterized by
tangent vectors $\eta \in T_x\mathcal M$ as
\begin{align}
  \gamma(t)=R_x(t\eta)\in \mathcal M, \text{such that~} (\gamma(0),\dot\gamma(0))=(x,\eta).
\end{align}
When $\mathcal M = \mathbb R^n$, the retraction is simply $R_x(\eta) = x + \eta.$
Other examples of retractions on the unit quaternion manifold $\mathcal Q$ and the unit sphere $S^2$ are as follows.
\begin{align}\label{eq:quaternion_retraction}
  R_q(v)= q\otimes \mathrm{Exp}(\phi),~
  \mathrm{Exp}(\phi)=
  \begin{bmatrix}
    \cos\|\phi\|\\
    \mathrm{sinc}(\|\phi\|)\,\phi
  \end{bmatrix},
\end{align}
\begin{align}\label{eq:sphere_retraction}
  R_s(w)= s\cos\|w\| + w~\mathrm{sinc}\norm{w}, 
\end{align}
where
$q \in \mathcal{Q}$, $\phi \in \mathbb{R}^3$, $s \in S^2$, $w \in T_s S^2$,
$v=q\otimes \phi \in T_q\mathcal Q$ and $\mathrm{sinc}(x)\coloneqq \sin x/x$. The symbol $\otimes$ denotes the quaternion multiplication.

\textbf{Vector transport:}
A vector transport $\mathcal T_{x\to y}:T_x\mathcal M\to T_y\mathcal M$ 
is used to compare tangent vectors based at different points.
Given a retraction $y=R_x(\eta)$, the retraction-induced transport is defined by
\begin{align}
  \mathcal T_{x\to y}(\delta\eta)\coloneqq DR_x(\eta)[\delta\eta]\in T_y\mathcal M,
\end{align}
which satisfies $\mathcal T_{x\to x}=\mathrm{id}$. 
This map serves as a computationally efficient approximation of parallel transport, avoiding the integration of differential equations.

\textbf{Covariant time derivative:}
For a curve $x(t)\in\mathcal M$ and a tangent vector field $\eta(t)\in T_{x(t)}\mathcal M$ along the curve,
$D\eta/dt$ represents the covariant derivative along $x(t)$, which gives an intrinsic notion of time variation
for vectors in moving tangent spaces.
Unless otherwise stated, $D\eta/dt$ denotes the covariant derivative
along $x(t)$ induced by the Levi--Civita connection of the Riemannian metric on $\mathcal M$.

\section{Problem Statement}
We address the following optimal control problem: 
\begin{subequations}\label{eq:optimal_control_problem}
\begin{align}
   \min_{x, u} \quad &J = \phi(x(t_f)) + \int_0^{t_f} L(x(t), u(t)) dt \label{eq:OCP_cost}\\
   \text{s.t.} \quad &\dot{x}(t) = f(x(t), u(t)), ~~ x(t) \in \mathcal{M}, ~~ u(t) \in \mathcal{U}, \label{eq:dynamics_constraint} \\
   &g(x(t), u(t)) \leq 0, \label{eq:path_constraint} \\
   &\psi(x(0), x(t_f)) = 0, \label{eq:boundary_condition}
\end{align}
\end{subequations}
where $\mathcal{M}$ and $\mathcal{U}$ are Riemannian manifolds.
Equations~\eqref{eq:dynamics_constraint}, \eqref{eq:path_constraint}, and \eqref{eq:boundary_condition} represent the system dynamics, path constraints, and boundary conditions, respectively.

The objective of this paper is to develop an \textit{intrinsic} numerical scheme for \eqref{eq:optimal_control_problem}, inspired by intrinsic SCvx,  \cite{kraisler_iscvx}. 
We combine a PS transcription with SCvx while respecting the geometry of $\mathcal{M}$ and $\mathcal{U}$.

\section{Classical Pseudospectral Method and Its Limitations}
\subsection{Standard Pseudospectral Method}\label{sec:standard_PS}
We first review a standard multi-segment ($hp$ method) PS transcription in Euclidean space.
PS methods are widely used in optimal control problems, 
owing to their properties as mentioned by \cite{Michael_Ross_review}: 
1) spectral convergence, and 2) mitigation of Runge's phenomenon.

In the $hp$ framework, the time horizon $[t_0, t_f]$ is partitioned into $N$ 
segments $[t_{h-1}, t_h]$ for $h=1, \dots, N$. The state and control trajectories are 
 approximated by a piecewise polynomial function. Within each segment $h$, 
we define a local time $\tau \in [-1, 1]$ and a set of flipped Radau nodes $\{\tau_i\}_{i=0}^p$ 
where $\tau_0 = -1$ following \cite{garg_RadauPS}. The corresponding nodal values are denoted 
as $\{x_i^h, u_i^h\}_{i=0}^p$. Using Lagrange polynomials $\mathcal{L}_k(\tau)$, the trajectories 
in each segment $h$ are approximated as:
\begin{align}\label{eq:approximation_hp}
x^h(\tau) \approx \sum_{0\leq k \leq p} x_k^h \mathcal{L}_k(\tau), \qquad
u^h(\tau) \approx \sum_{0\leq k \leq p} u_k^h \mathcal{L}_k(\tau).
\end{align}

It is important to note that while the state trajectory is 
enforced to be continuous across segments, the control trajectory is generally discontinuous 
at the segment boundaries. To ensure state continuity, the following linking conditions 
are imposed for $h=1, \dots, N-1$:
\begin{align}
  \label{eq:linking_condition}x_p^h = x_0^{h+1}.
\end{align}
Under this $hp$ discretization, 
the collocation dynamics for each segment $j$ at collocations $i=1, \dots, p$ are 
given by:
\begin{align}\label{eq:collocation_constraint_hp}
  \sum_{0\leq j \leq p} D_{ij}x_j^h = \sigma^h f(x_i^h, u_i^h), 
  \qquad \sigma^h \coloneqq \frac{t_h - t_{h-1}}{2},
\end{align}
where $D_{ij} = \mathcal{L}_j'(\tau_i)$ is the differentiation matrix. 
In this work, we assume uniform segmentation, 
yielding a common time-scaling factor for every segment $\sigma^h = \sigma = (t_f - t_0)/(2N)$.

\subsection{Pseudospectral Successive Convexification}
The PS transcription yields a finite-dimensional nonlinear program (NLP) 
with decision variables
$\vb*{x}\coloneqq\{x_i^{h}\}_{i \in [1, p]}^{h\in [1,N]}$, $\vb*{u}\coloneqq\{u_i^{h}\}_{i \in [1, p]}^{h\in [1,N]}$,
and optionally $\sigma$.
Using the quadrature weights $w_i$, the integral cost \eqref{eq:OCP_cost} is expressed by
\begin{align}
  J
  &= \phi\!\left(x_{p}^{N}\right)
   + \sigma\!\! \sum_{ 1 \leq h \leq N} \sum_{{ 1 \leq i \leq p}}
     w_i\, L\!\left(x_i^{h},u_i^{h}\right).
\end{align}

Other constraints \eqref{eq:dynamics_constraint}--\eqref{eq:boundary_condition} and linking condition are transcribed as
\begin{subequations}\label{eq:hp_nlp_constraints}
\begin{align}
  \sum_{0\leq j \leq p}\! D_{ij}^{h} x_j^{h}
  &= \sigma\, f\!\left(x_i^{h},u_i^{h}\right),
  && \substack{i=1,\dots,p\\ h=1,\dots,N} \label{eq:hp_collocation}\\
  g\!\left(x_i^{h},u_i^{h}\right)
  &\le 0,
  && \substack{i=1,\dots,p\\ h=1,\dots,N} \\
  \psi\!\left(x_0^{1},x_{p}^{N}\right)
  &= 0, & & \\
  x_{p}^{h}
  &= x_{0}^{h+1}, && {\scriptstyle h=1,\dots,N-1}. \label{eq:hp_linking}
\end{align}
\end{subequations}

This NLP is generally nonconvex due to the nonlinear dynamics and constraints.
\cite{sagliano_SPsCvx} proposes an SCvx approach to solve this NLP efficiently.
At each iteration, it linearizes nonconvex constraints about a reference
$\bigl(\bar{x}_i^{h}, \bar{u}_i^{h}, \bar{\sigma}\bigr)$.
Let $\Delta\sigma \coloneqq \sigma-\bar{\sigma}$, and let $A_i^{h}$ and $B_i^{h}$ denote the Jacobians of $f$
w.r.t. $x$ and $u$ at $(\bar{x}_i^{h},\bar{u}_i^{h})$.
Then, for all collocation points $\substack{i=1,\dots,p\\ h=1,\dots,N}$,
\begin{align}\label{eq:hp_linearized_dynamics_rhs}
  \begin{split}
    \sum_{0\leq j \leq p} \!\!\!D_{ij} x_j^{h}
= \bar{\sigma}\bigl( f_i^{h} + A_i^{h}\Delta x_i^{h} + B_i^{h}\Delta u_i^{h}\bigr)+ f_i^{h}\,\Delta\sigma,
  \end{split}
\end{align}
where $f_i^{h} \coloneqq f(\bar{x}_i^{h},\bar{u}_i^{h})$,
$\Delta x_i^{h} \coloneqq x_i^{h} - \bar{x}_i^{h}$, and 
$\Delta u_i^{h} \coloneqq u_i^{h} - \bar{u}_i^{h}$.
The path constraints are linearized similarly, and the resulting convex subproblem
is solved iteratively.

\subsection{Geometric Inconsistency on Manifolds}

Classical transcription in Euclidean space $\mathbb{R}^n$ becomes 
geometrically inconsistent on nonlinear manifolds.
The main issues of \eqref{eq:hp_linearized_dynamics_rhs} are as follows:
\begin{itemize}
  \item \textbf{Violation of state feasibility:}
  Ambient-space sums $\sum_{k} x_k \mathcal{L}_k(\tau)$ do not guarantee the interpolant lies on $\mathcal{M}$, even if $x_k\in\mathcal{M}$.

  \item \textbf{Vector space mismatch:}
  $\sum_{j} D_{ij}x_j$ is an element of $\mathbb{R}^n$,
  whereas the vector field satisfies $f(x_i,u_i)\in T_{x_i}\mathcal{M}$.
  These objects live in different spaces unless one explicitly projects or transports them.

  \item \textbf{Ill-posed increments for linearization:}

  The Euclidean increment $\Delta x \!=\! x\!-\!\bar{x}$ is not intrinsic to $\mathcal{M}$, and ambient-space linearizations (e.g., $f(x,u)\! \approx \!f(\bar{x},\bar{u})\! +\! A\Delta x \!+\! B\Delta u$) do not necessarily lie in $T_{x}\mathcal{M}$.
\end{itemize}
Existing approaches, such as \cite{sagliano_accD}, attempted to guarantee geometric consistency extrinsically, 
by imposing additional constraints or introducing auxiliary variables.
In contrast, this paper develops an \textit{intrinsic} scheme that strictly respects the geometry of $\mathcal{M}$ and $\mathcal{U}$.

\section{Intrinsic Pseudospectral Successive Convexification}
\begin{figure}[bt]
    \centering
    \includegraphics[width=0.95\linewidth]{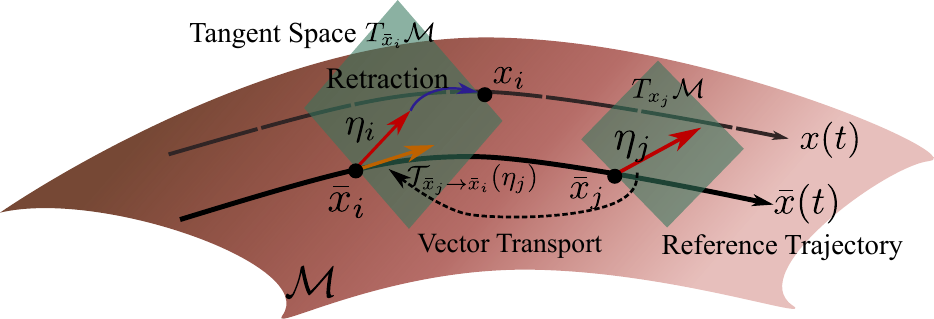} 
    \caption{Concept of iPS SCvx method. The trajectories are parameterized by retraction maps, and the perturbation dynamics are defined on the tangent bundle.}
    \label{fig:concept_iPS}
\end{figure}

To resolve the geometric inconsistencies, the problem is formulated on the tangent bundle.
Since tangent spaces admit vector-space operations, classical PS methods become directly applicable.
Parameterizing the state as a perturbation from a reference trajectory and applying convexification yield the proposed iPS SCvx scheme.
Figure~\ref{fig:concept_iPS} illustrates the overall concept.

\subsection{Perturbation Dynamics on Tangent Bundle}
Near a reference $(\bar{x}(t),\bar{u}(t)) \in \mathcal{M}\times\mathcal{U}$, trajectories $(x(t), u(t))$ are parameterized via retractions $R$ as:
\begin{subequations}\label{eq:retraction_parameterization}
\begin{align}
  x(t) &= R_{\bar{x}(t)}(\eta(t)), \qquad \eta(t)\in T_{\bar{x}(t)}\mathcal{M},\\
  u(t) &= R_{\bar{u}(t)}(\xi(t)), \qquad \xi(t)\in T_{\bar{u}(t)}\mathcal{U}.
\end{align}
\end{subequations}
Since $\eta(t)$ resides in the time-varying tangent space $T_{\bar{x}(t)}\mathcal{M}$, its rate of change is characterized by the covariant derivative $D\eta/dt$.
Time-differentiation of \eqref{eq:retraction_parameterization} yields the velocity expansion:
\begin{align}\label{eq:x_dot_expansion}
  \dot{x}(t)
  = D_1R_{\bar{x}(t)}(\eta(t))[\dot{\bar{x}}(t)]
  + D_2R_{\bar{x}(t)}(\eta(t))\!\left[\frac{D\eta}{dt}\right],
\end{align}
where $D_1R$ and $D_2R$ denote the partial derivatives of $R_x(\eta)$ with respect to the base point $x$
and the tangent vector $\eta$, respectively.
Substituting $\dot{x}(t)=f(x(t),u(t))$ into \eqref{eq:x_dot_expansion} gives the perturbation dynamics on $T\mathcal{M}$:
\begin{align}\label{eq:perturbed_dynamics_retraction}
\begin{split}
    \frac{D\eta}{dt} &= \mathcal F(\eta,\xi)\\
&\! \coloneqq \!\bigl(D_2R_{\bar{x}}(\eta)\bigr)^{\raisebox{0.3ex}{\kern -0.3em $\scriptstyle -1$}}\!\Bigl[f\!\bigl(R_{\bar{x}}(\eta),R_{\bar{u}}(\xi)\bigr)
      \!- \!D_1R_{\bar{x}}(\eta)[\dot{\bar{x}}]\Bigr].
\end{split}
\end{align}
Here $D_2R_{\bar{x}}(\eta)$ is assumed to be invertible in a neighborhood of $0_{\bar{x}}$.
This holds locally since
$D_2R_{\bar{x}}(0)=\mathrm{Id}$ by \eqref{eq:retraction_def}.

\subsection{Discretization and Linearization of Dynamics and Constraints}\label{sec:disc_lin}
Equation~\eqref{eq:perturbed_dynamics_retraction} is transcribed by an $hp$-pseudospectral scheme and linearized analogously to Section~\ref{sec:standard_PS},
 focusing on pertubation vectors $\eta$ and $\xi$ which live in the tangent spaces.
Since $\eta_k^{h}\in T_{\bar{x}_k^{h}}\mathcal{M}$ live in different tangent spaces, interpolation is formed after transporting nodal perturbations by the retraction-induced vector transport $\mathcal{T}$.
On each segment $h$, $\eta^{h}(\tau)$ is approximated by
\begin{equation}\label{eq:eta_interpolant}
  \eta^{h}(\tau)
  \approx
  \sum_{0 \leq k \leq p}
 \qty[ \mathcal{T}_{\bar{x}_k^{h}\to \bar{x}^{h}(\tau)}\!\bigl(\eta_k^{h}\bigr)\,\mathcal{L}_k(\tau)].
\end{equation}

Taking the covariant derivative at $t=t_i$ yields
\begin{equation}\label{eq:covariant_diff_discrete}
  \left.\frac{D\eta^{h}}{dt}\right|_{t_i}
  =
  \frac{1}{\sigma}\sum_{0 \leq k \leq p} \mqty[D_{ik}\,
  \mathcal{T}_{\bar{x}_k^{h}\to \bar{x}_i^{h}}\!\bigl(\eta_k^{h}\bigr)]
  + \mathcal{E}^{h}_{i}\!\bigl(\eta_i^{h}\bigr),
\end{equation}
where $\mathcal{E}^{h}_{i}$ is the linear correction for vector transport:
\begin{equation}\label{eq:E_transport_def}
  \mathcal{E}^{h}_{i}\!\bigl(\eta\bigr)
  \coloneqq
  \left.\frac{D}{dt}\,
  \mathcal{T}_{\bar{x}_i^{h}\to \bar{x}^{h}(t)}(\eta)
  \right|_{t_i}.
\end{equation}

Substituting \eqref{eq:covariant_diff_discrete} into \eqref{eq:perturbed_dynamics_retraction}
gives the equations on $T_{\bar{x}_i^{h}}\mathcal{M}$:
\begin{align}\label{eq:eta_collocation}
  \sum_{0 \leq k \leq p}\!\!\! \qty[D_{ik}\,
  \mathcal{T}_{\bar{x}_k^{h}\to \bar{x}_i^{h}}\!\bigl(\eta_k^{h}\bigr)]
  =
  \sigma\qty(\mathcal{F}\!\left(\eta_i^{h},\xi_i^{h}\right) - \mathcal{E}^{h}_{i}\!\bigl(\eta_i^{h}\bigr)).
\end{align}

Note that the left-hand side of \eqref{eq:eta_collocation} is a linear map with respect to $\{\eta_i^{h}\}_{i=0}^p$.
Then the right-hand side is linearized about $(\eta,\xi, \sigma)=(0,0, \bar{\sigma})$, which gives

\begin{align}\label{eq:linearized_collocation}
      \sum_{0 \leq k \leq p}\!\!\!\! \qty[D_{ik}\!
  \mathcal{T}_{\bar{x}_k^{h}\!\to \bar{x}_i^{h}}\!\bigl(\eta_k^{h}\bigr)]
  \!\!=\!
  \bar{\sigma}\qty[\!\rho_i^{h} \!\!+\!\!
   \tilde{A}_i^{h}\!\bigl(\eta_i^{h}\bigr)\!\! + \!\!B_i^{h}\!\bigl(\xi_i^{h}\bigr)\! ]\!\!+\! \rho_i^{h}\Delta\sigma.
\end{align}
Here, the operators $\tilde{A}_i^{h}$, $B_i^{h}$ are defined by
\begin{subequations}\label{eq:ABCS_defs}
  \begin{align}
  \tilde{A}_i^{h}(\cdot)
  &\coloneqq
  \Bigl(\mathrm{D}_x f(\bar{x}_i^{h},\bar{u}_i^{h}) - C_i^{h} + S_{i, \rho_i^{h}}^{h} - \mathcal{E}_i^h\Bigr)(\cdot),\\
   B_i^{h}(\cdot)
  &\coloneqq \mathrm{D}_u f(\bar{x}_i^{h},\bar{u}_i^{h})(\cdot),\\
  C_i^{h}(\cdot)
  &\coloneqq
  \left.
  \mathrm{D}_{\eta}\Bigl(D_1R_{\bar{x}_i^{h}}(\eta)[\dot{\bar{x}}_i^{h}]\Bigr)
  \right|_{\eta=0}(\cdot),\\
  S_{i, \rho_i^{h}}^{h}(\cdot)
  &\coloneqq
  \biggl(\left.\mathrm{D}_{\eta}\Bigl(\bigl(D_2R_{\bar{x}_i^{h}}(\eta)\bigr)^{\raisebox{0.3ex}{\kern -0.2em $\scriptstyle -1$}}\Bigr)
  \right|_{\eta=0}(\cdot)\biggr)(\rho_i^{h}).
\end{align}
\end{subequations}
In \eqref{eq:ABCS_defs}, $\mathrm{D}_x$, $\mathrm{D}_u$, and $\mathrm{D}_\eta$ denote partial differentials with respect to the state, control, and tangent perturbation, respectively.
The term $\rho_i^{h}$ represents the reference defect 
and reference velocity $\dot{\bar{x}}_i^{h}$ is given by
\begin{align}\label{eq:defect_rho}
    \rho_i^{h}
  \coloneqq
  f\!\left(\bar{x}_i^{h},\bar{u}_i^{h}\right)
  - \dot{\bar{x}}_i^{h}, ~
    \dot{\bar{x}}_i^{h} \approx \frac{1}{\bar{\sigma}}\sum_{0 \leq k \leq p} D_{ik}\,
    R_{\bar{x}_i^{h}}^{-1}\!\left(\bar{x}_k^{h}\right),
\end{align}
where $R_{\bar{x}_i^{h}}^{-1}(\cdot)$ is the inverse retraction map.

Nonconvex path constraints $g(x,u)\le 0$ are also linearized using retraction as
\begin{align}\label{eq:path_linearization}
  g\!\left(\bar{x}_i^{h},\bar{u}_i^{h}\right)
  + G^{x,h}_i\!\bigl(\eta_i^{h}\bigr)
  + G^{u,h}_i\!\bigl(\xi_i^{h}\bigr)
  \le 0.
\end{align}
where $G_{i}^{x,h} \coloneqq \mathrm{D}_x g(\bar{x}_i^{h},\bar{u}_i^{h})$ and $G_{i}^{u,h} \coloneqq \mathrm{D}_u g(\bar{x}_i^{h},\bar{u}_i^{h})$.

\subsection{$hp$-method from a Geometric Perspective}\label{sec:hp_geometric}
As discussed in Section~\ref{sec:standard_PS}, the $hp$-method divides the time horizon into multiple segments and connects states by linking conditions.
From a geometric viewpoint, the $hp$-method promotes locality.
Within each short segment, the reference nodes are closer, so the retraction-induced transports
$\mathcal T_{\bar x_k^h\to \bar x_i^h}$ tend to remain near the identity and the discrete operator
$\sum_k D_{ik}\mathcal T_{\bar x_k^h\to \bar x_i^h}(\eta_k^h)$ is less sensitive to curvature.
Segmentation also helps keep iterates within a neighborhood where $R_{\bar x}$ and $D_2R_{\bar x}$ are well behaved.

State continuity \eqref{eq:hp_linking} at segment interfaces is enforced through the retraction parameterization:
\begin{align}\label{eq:linking_condition_retraction}
  R_{\bar{x}_p^{h}}\!\left(\eta_p^{h}\right)
  =
  R_{\bar{x}_0^{h+1}}\!\left(\eta_0^{h+1}\right).
\end{align}
The reference trajectory is initialized and updated to satisfy continuity at interfaces (i.e., $\bar{x}_p^{h}=\bar{x}_0^{h+1}$).
Under this condition, both perturbations $\eta_p^{h}$ and $\eta_0^{h+1}$ reside in the same tangent space $T_{\bar{x}_p^{h}}\mathcal{M}$.
Since the retraction $R_x(\cdot)$ is locally diffeomorphic at the origin, \eqref{eq:linking_condition_retraction} simplifies to
\begin{align}\label{eq:eta_link}
  \eta_p^{h} = \eta_0^{h+1}.
\end{align}
This linear formulation avoids the need for nonlinear equality constraints at the segment boundaries, significantly improving numerical tractability.

\subsection{Coordinate Representation via Local Frames}\label{sec:coords}
To solve the convex subproblem numerically, 
abstract tangent-space quantities are converted into concrete coordinate representations.
Following \cite{kraisler_iscvx}, orthonormal frames $E_i^{h}$ and $F_i^{h}$ are introduced for the tangent spaces $T_{\bar{x}_i^{h}}\mathcal{M}$ and $T_{\bar{u}_i^{h}}\mathcal{U}$, respectively.
These frames map the perturbations $\eta_i^h, \xi_i^h$ to coordinate vectors $\hat{\eta}_i^{h}\in\mathbb{R}^n, \hat{\xi}_i^{h}\in\mathbb{R}^m$ such that  $\eta_i^h = E_i^h\hat\eta_i^h$ and $\xi_i^h = F_i^h\hat\xi_i^h$.

The operators and transport maps are converted into matrix forms $[\tilde{A}]_i^{h}, [{B}]_i^{h}, [{G}^x]_i^{h}, [{G}^u]_i^{h}, [{T}]_{ik}^{h}$ via projection onto these frames.
Similarly, $\hat{\rho}_i^{h}$ denote the coordinate representations of the defect vector.

Consequently, the discretized constraints \eqref{eq:linearized_collocation} and \eqref{eq:path_linearization} are expressed in coordinates as:
\begin{align}\label{eq:linear_coord}
  \begin{split}
    \sum_{0\leq k \leq p} D_{ik}\,[{T}]_{ik}^{h}\,\hat{\eta}_k^{h}
  &= \bar{\sigma}\bigl(\hat{\rho}_i^{h} + [\tilde{A}]_i^{h}\hat{\eta}_i^{h}+[{B}]_i^{h}\hat{\xi}_i^{h}\bigr)\\[-2ex]
  &\quad  + \hat{\rho}_i^{h}\Delta\sigma + \hat{\nu}_i^{h}, 
  \end{split}\\
   g(\bar{z}_i^{h}) &+ [{G}^x]^{h}_i\,\hat{\eta}_i^{h} + [{G}^u]^{h}_i\,\hat{\xi}_i^{h} \le \hat{s}_i^{h},
\end{align}
where $\hat{\nu}_i^{h}$ and $\hat{s}_i^h$ are virtual control and slack variables added to ensure feasibility as in \cite{kraisler_iscvx}.

Finally, to simplify the linking conditions, a consistent frame orientation is enforced at segment interfaces (i.e., $E_p^{h}=E_0^{h+1}$).
This reduces the geometric continuity condition to a simple linear equality in $\mathbb{R}^n$:
\begin{align}
  \hat{\eta}_p^{h}=\hat{\eta}_0^{h+1}.
\end{align}

\subsection{Formulation of the Convex Subproblem}
At iteration $\ell$, we optimize the perturbations $\mathcal{Z} \coloneqq \{\hat{\eta}, \hat{\xi}, \Delta\sigma, \hat{\nu}, \hat{s}, r\}$ by solving the following subproblem:
\begin{subequations} \label{eq:convex_subproblem}
\begin{align}
\min\limits_{\mathcal{Z}}&\quad
J_{\mathrm{cvx}}^{(\ell)}
+ \sum\limits_{\mathclap{\substack{h=1,\dots,N\\ i=1,\dots,p}}} \! w_i\Bigl(
\mu_\nu \|\hat{\nu}_i^{h}\|_1 + \mu_s(\hat{s}_i^{h})_{+} + \mu_r r_i^{h}
\Bigr)\\
\begin{split}
  \text{s.t.} \quad & \sum_{0 \leq k \leq p} \!\!\! D_{ik}\,[{T}]_{ik}^{h}\,\hat{\eta}_k^{h}=
\bar{\sigma}\bigl(\hat{\rho}_i^{h}  + [\tilde{A}]_i^{h}\hat{\eta}_i^{h} +[{B}]_i^{h}\hat{\xi}_i^{h}\bigr)\\[-2ex]
&\qquad\qquad\qquad + \hat{\rho}_i^{h}\Delta\sigma
  + \hat{\nu}_i^{h}\label{eq:sub_dyn}
\end{split}\\
    &g(\bar{x}_i^{h}, \bar{u}_i^{h}) + [{G}^x]^{h}_i\,\hat{\eta}_i^{h} + [{G}^u]^{h}_i\,\hat{\xi}_i^{h} \le \hat{s}_i^{h}, \label{eq:sub_path}\\
    &  \hat{\eta}_p^{h}=\hat{\eta}_0^{h+1}, \label{eq:sub_link}\\
    & \norm{\hat{\xi}_i^{h}}_{2}^2 \le r_i^{h}, \label{eq:trust_region}\\
    &  \psi\Bigl(\hat{\eta}_0^{1}, \hat{\eta}_p^{N}\Bigr)= 0
\end{align}
\end{subequations}

Here, $J_{\mathrm{cvx}}^{(\ell)}$ 
is assumed to be a convex quadratic approximation of the trajectory cost about the reference (see \cite{kraisler_iscvx}).
The weights $\mu_{\nu}, \mu_{s}, \mu_r > 0$ penalize the virtual controls, constraint violations, and the trust region radius, respectively.

Constraint~\eqref{eq:trust_region} enforces a trust region on the control step.
Since $\mathcal U$ need not provide explicit bounds, the convexified subproblem may propose
$\hat\xi_i^h$ that leaves the neighborhood where the retraction $R_{\bar u_i^h}(\xi)$ and its
linearization are reliable. The step is therefore bounded by $\|\hat\xi_i^h\|_2^2 \le r_i^h$ and
penalizes $r_i^h$ in the objective to discourage overly large updates.
Since $\xi_i^h$ is expressed in an orthonormal frame $F_i^h$ of $T_{\bar u_i^h}\mathcal U$,
the coordinate norm equals the intrinsic Riemannian norm:
$\ \|\xi_i^h\|_{\bar u_i^h}=\|\hat\xi_i^h\|_2$.
With a non-orthonormal basis, $\|\xi\|^2_{\bar{u}}=\hat\xi^\top G\hat\xi$, $G_{ab}=\langle F_a,F_b\rangle$.

The overall procedure is summarized in Algorithm \ref{alg:iscvx}.
\begin{algorithm}
\caption{Intrinsic Pseudospectral Successive Convexification}
\label{alg:iscvx}
\begin{algorithmic}[1]
\Require Initial trajectory $(\bar{x}^0, \bar{u}^0, \bar{\sigma}^0) \in \mathcal{M}\times \mathcal{U} \times \mathbb{R}$, 
Weights $\mu_{\nu}, \mu_{s}, \mu_r > 0$, convergence tolerance $\varepsilon > 0$.
\State $k \gets 0$
\Repeat
    \State \textbf{Matrices Calculation:}
    \Statex \qquad Compute $[\tilde{A}], [{B}]$, $[G^x]$, $[G^u]$, and $[T]$ from
    \Statex \qquad  reference trajectory $(\bar{x}^k, \bar{u}^k, \bar{\sigma}^k)$.
    
    \State \textbf{Convex Optimization:}
    \State \quad Solve the convex subproblem : \eqref{eq:convex_subproblem}
    \State \quad Let $(\hat{\eta}^*, \hat{\xi}^*, \Delta \sigma^*)$ be the optimal solution.

    \State \textbf{Update:}
    \State \quad $\bar{x}_i^{k+1} \gets R_{\bar{x}_i^{k}}(E_i \hat{\eta}_i^*)$, $\bar{u}_i^{k+1} \gets R_{\bar{u}_i^{k}}(F_i \hat{\xi}_i^*), \quad \forall i$
    \State \quad $\bar{\sigma}^{k+1} \gets \bar{\sigma}^{k} + \Delta \sigma^*$
    
    \State $k \gets k + 1$
\Until{$\|\hat{\eta}^*\| < \varepsilon$}
\Statex \hspace*{-\algorithmicindent}\Return $(\bar{x}^*, \bar{u}^*, \bar{\sigma}^*)$
\end{algorithmic}
\end{algorithm}

\section{Numerical Example: Six Degree of Freedom Landing Guidance}
This section demonstrates the proposed iPS SCvx method on a six-degree-of-freedom rocket landing problem.

\subsection{Problem Formulation}
Following \cite{sagliano_accD}, we fix $t_f=\SI{4}{Ut}$ and neglect RCS-induced mass depletion and torques. 
The formulation is briefly summarized below:
\begin{subequations}
\begin{align}
  & \hspace{-10mm}\min_{T_\mathrm{mag}, \vb*{u}_{\mathrm{dir}}} \quad  J = -m(t_f) \label{eq:cost} \\
  \text{s.t.} \quad & \left. \begin{aligned}
    & \dot{m} = -\alpha T_\mathrm{mag}, \quad \dot{\vb*{r}} = \vb*{v} \\
    & \dot{{v}} = \frac{T_\mathrm{mag}}{m}C_{\mathcal{I}\mathcal{B}}{u}_{\mathrm{dir}} + {g} + \frac{{D}}{m},\;\dot{{q}} = \frac{1}{2} {q}\otimes {\omega} \\
    & \dot{{\omega}} = J^{-1}\left( {\ell}_\mathrm{arm}\times {T_\mathrm{mag}} {u}_{\mathrm{dir}} - {\omega} \times J {\omega}\right) 
  \end{aligned} \right\} \\
  & \left. \begin{aligned}
    & \sqrt{r_y^2 + r_z^2} \leq r_x \cot \gamma, \quad \norm{{\omega}} \leq \omega_\mathrm{max} \\
    & T_\mathrm{min} \leq T_\mathrm{mag} \leq T_\mathrm{max}, \quad m(t_f) \geq m_\mathrm{dry}
  \end{aligned} \right\} \\
  & u_{\mathrm{dir},1} \geq \cos\delta_{\mathrm{max}}, \quad q_y^2 + q_z^2 \leq \sin^2 \qty({\phi_\mathrm{max}}/{2}) \label{eq:attitude_constraint}
\end{align}
\end{subequations}

The state vector ${x} = [m, {r}^\top, {v}^\top, {q}^\top, {\omega}^\top]^\top$ 
includes mass, position, velocity, quaternion, and angular velocity. 
The control input ${u} = [T_{\mathrm{mag}}, {u}_{\mathrm{dir}}^\top]^\top$ 
consists of the thrust magnitude and direction vector. 
For the manifold components $q \in \mathcal{Q}$ and $u_{\mathrm{dir}} \in S^2$, 
the retractions \eqref{eq:quaternion_retraction} and \eqref{eq:sphere_retraction} are applied. 
Constraints on these variables \eqref{eq:attitude_constraint} are linearized via 
the intrinsic scheme in Section \ref{sec:disc_lin}, 
with a trust region \eqref{eq:trust_region} imposed on $u_{\mathrm{dir}}$ updates.

All boundary conditions and constants follow \cite{sagliano_accD}, except that we set the initial attitude to
$\vb*{q}(0)=[0.7428,-0.04278,0.03559,0.6672]^\top$ to enable a fair comparison with an extrinsic ACCD baseline
under the same boundary/parameter settings.
The penalty weights are $\mu_\nu=\SI{e4}{},\ \mu_s=\SI{1e-1}{},\ \mu_r=\SI{3e-2}{}$.
We use an $hp$ scheme with $N=5$ segments and $p=10$ flipped Radau collocations per segment.

Simulations were run on a 1.4 GHz Intel i5 PC (8GB RAM) using CVX/ECOS in MATLAB 2025b.

\subsection{Simulation Results}

\begin{figure}[bt]
  \centering
  \begin{subfigure}[t]{0.99\linewidth}
    \centering
    \includegraphics[width=0.75\linewidth]{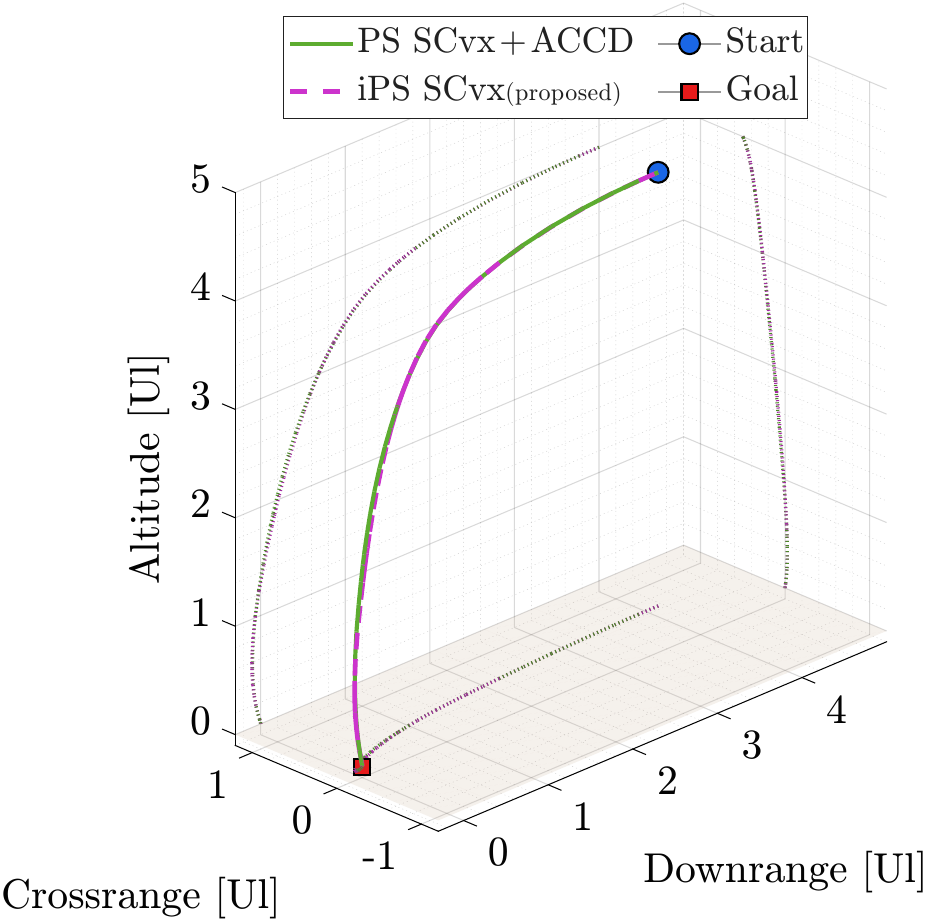} 
    \caption{Trajectories}
    \label{fig:trajectory}
  \end{subfigure}
    \begin{subfigure}[t]{0.492\linewidth}
    \centering
    \includegraphics[width=\linewidth]{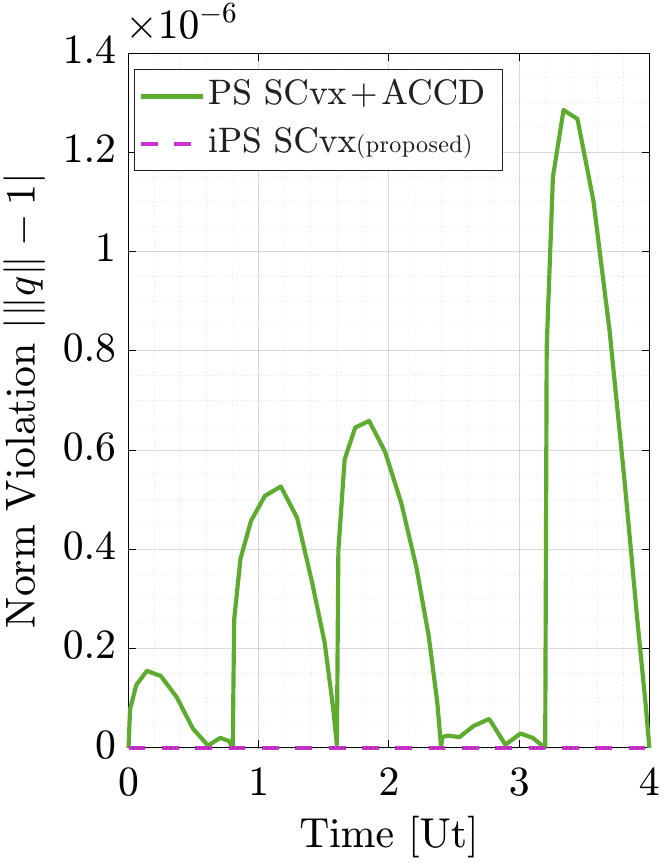} 
    \caption{Violation from $\norm{{q}} = 1$}
    \label{fig:q_norm}
  \end{subfigure}
    \begin{subfigure}[t]{0.492\linewidth}
    \centering
    \includegraphics[width=\linewidth]{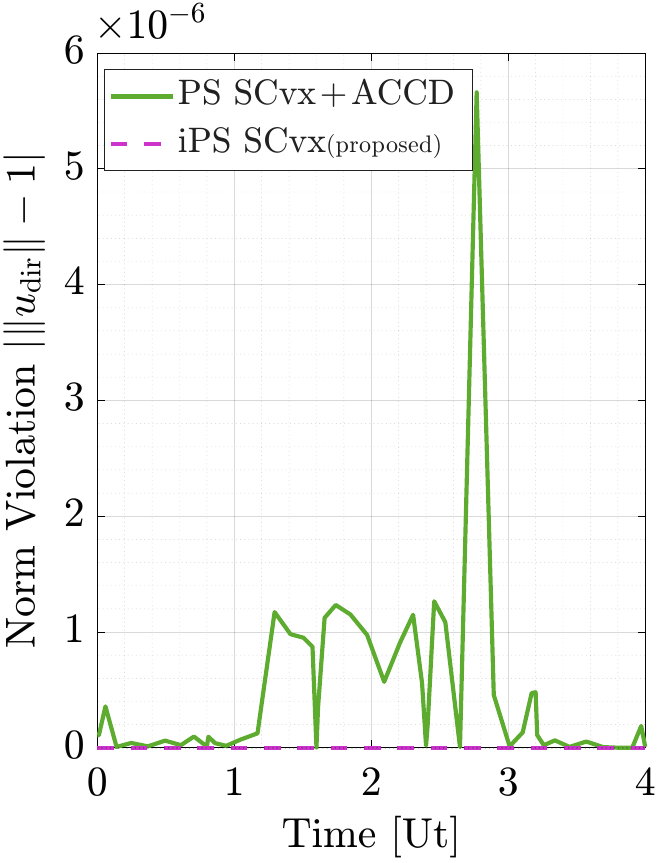} 
    \caption{Violation from $\norm{{u}_\mathrm{dir}} = 1$}
    \label{fig:u_dir_norm}
  \end{subfigure}
 \caption{Comparison of PS SCvx (ACCD) vs. proposed iPS SCvx. (a) Trajectories. (b)(c) Unit-norm violations for quaternion $\|{q}\|$ and thrust direction $\|{u}_\mathrm{dir}\|$ on a linear scale. The proposed method stays below $10^{-15}$.}
  \label{fig:results_comparison}
\end{figure}

Figure~\ref{fig:results_comparison} compares the PS SCvx method with extrinsic ACCD and the proposed iPS SCvx method.
As shown in Fig.~\ref{fig:trajectory}, both methods converge to essentially the same solution in terms of trajectory shape and terminal mass
($m_f = \SI{1.954}{kg}$ for PS SCvx with ACCD and $m_f = \SI{1.953}{kg}$ for iPS SCvx), indicating almost identical optimality.
In contrast, Figs.~\ref{fig:q_norm}--\ref{fig:u_dir_norm} highlight a key difference in geometric consistency.
The PS SCvx with ACCD method exhibits unit-norm drift (up to $10^{-6}$), whereas the iPS SCvx method preserves $\|{q}\|=1$ and $\|{u}_\mathrm{dir}\|=1$ by construction, 
remaining at the level of machine precision ($\sim10^{-15}$).
This result is gained without renormalization or additional manifold-embedding constraints.

As shown in Table \ref{tab:computation_time}, 
the proposed iPS SCvx method demonstrates superior computational 
performance compared to the PS SCvx method with ACCD. 
It requires fewer iterations and achieves a reduction in total computation time, 
alongside a lower per-iteration time. 
The improved convergence rate is likely due to the intrinsic formulation's better 
handling of manifold geometry, which leads to more accurate linearizations and faster progress 
towards the optimal solution. 
Furthermore, the reduction in per-iteration time is 
attributed to a decrease in the number of decision variables and constraints. 
Although iPS SCvx introduces additional matrix computations, 
the computational savings from the reduced problem dimensionality
outweigh this overhead.
\begin{table}[bt]
  \centering
  \caption{Computational performance of PS SCvx (ACCD) vs. proposed iPS SCvx. Numbers in parentheses denote standard deviation.}
  \label{tab:computation_time}
\begin{tabular}{@{}lccc@{}}
\toprule
          & Iter. & Total [s] &  Time$/$Iter. [s] \\ \midrule
PS SCvx + ACCD & 88        & \SI{1.341e3}{} (33.85)       & 15.23  (0.6022)          \\
iPS SCvx  & 16         &  \SI{1.680e2}{}  (29.81)       & 10.50      (2.281)          \\ \bottomrule
\end{tabular}
\end{table}

\section{Conclusion}
This paper introduced an iPS SCvx method, which bridges
pseudospectral transcription and intrinsic SCvx for manifold-constrained optimal control.
By expressing dynamics and increments intrinsically in local tangent spaces and coupling them with
retraction-induced transports, the proposed transcription preserves manifold feasibility by construction
and removes the need for extrinsic manifold constraints or auxiliary variables typically used in classical formulations.
A six-degree-of-freedom powered-landing example with unit-quaternions and 
thrust direction vectors demonstrated comparable optimality, alongside 
improved feasibility and computational efficiency.

The present study focuses on manifolds with
closed-form retractions and associated differential operators.
An important next step is to extend the implementation to general
manifolds by numerically evaluating the required retraction
differentials and transport terms, enabling its application without closed-form operators.
Moreover, future work 
will explore optimization algorithm to further reduce computation time for real-time applications.

\section*{DECLARATION OF GENERATIVE AI AND AI-ASSISTED TECHNOLOGIES IN THE WRITING PROCESS}
During the preparation of this work the authors used chatGPT in order to assist with English proofreading and code generation assistance.
After using this service, the authors reviewed and edited the content as needed and take full responsibility for the content of the publication.
\bibliography{ifacconf}             

@article{wojciech_stabilization,
title = {Methods for constraint violation suppression in the numerical simulation of constrained multibody systems – A comparative study},
journal = {Computer Methods in Applied Mechanics and Engineering},
volume = {200},
number = {13},
pages = {1568-1576},
year = {2011},
author = {Wojciech Blajer}
}

@book{Absil_optManifold,
 ISBN = {9780691132983},
 author = {P.-A. Absil and R. Mahony and R. Sepulchre},
 publisher = {Princeton University Press},
 address   = {Princeton, NJ},
 title = {Optimization Algorithms on Matrix Manifolds},
 urldate = {2025-12-13},
 year = {2008}
}

@article{kraisler_iscvx,
  author       = {Kraisler, Spencer and Mesbahi, Mehran and A{\c c}{\i}kme{\c s}e, Beh{\c c}et},
  title        = {Intrinsic Successive Convexification: Trajectory Optimization on Smooth Manifolds},
  journal = {IEEE Control Systems Letters},
  volume       = {9},
  pages        = {408--413},
  year         = {2025},
}

@conference{sagliano_SPsCvx,
author = {Sagliano, Marco and Heidecker, Ansgar and Hern{\'a}ndez, Jos{\'e} Mac{\'e}s and Far{\`i}, Stefano and Schlotterer, Markus and Woicke, Svenja and Seelbinder, David and Dumont, Etienne},
title = {Onboard Guidance for Reusable Rockets: Aerodynamic Descent and Powered Landing},
booktitle = {AIAA Scitech 2021 Forum},
year = {2021},
chapter = {},
pages = {}
}

@article{lu_CCD,
author = {Lu, Ping},
title = {Convex–Concave Decomposition of Nonlinear Equality Constraints in Optimal Control},
journal = {Journal of Guidance, Control, and Dynamics},
volume = {44},
number = {1},
pages = {4-14},
year = {2021},

}

@article{sagliano_accD,
author = {Sagliano, Marco and Seelbinder, David and Theil, Stephan and Lu, Ping},
year = {2024},
month = {01},
pages = {20-35},
title = {Six-Degree-of-Freedom Rocket Landing Optimization via Augmented Convex–Concave Decomposition},
volume = {47},
number = {1},
journal = {Journal of Guidance, Control, and Dynamics},
}

@conference{szmuk_Free6DoF,
author = {Szmuk, Michael and A\c{c}\i{}kme\c{s}e, Beh\c{c}et},
title = {Successive Convexification for 6-DoF Mars Rocket Powered Landing with Free-Final-Time},
booktitle = {2018 AIAA Guidance, Navigation, and Control Conference},
chapter = {},
pages = {},
year = {2018}
}

@ARTICLE{Saccon_LieGroupProjection,
  author={Saccon, Alessandro and Hauser, John and Aguiar, A. Pedro},
  journal={IEEE Transactions on Automatic Control}, 
  title={Optimal Control on Lie Groups: The Projection Operator Approach}, 
  year={2013},
  volume={58},
  number={9},
  pages={2230-2245},
 }

@ARTICLE{Bordalba_direct_collocation,
  author={Bordalba, Ricard and Schoels, Tobias and Ros, Lluís and Porta, Josep M. and Diehl, Moritz},
  journal={IEEE Transactions on Robotics}, 
  title={Direct Collocation Methods for Trajectory Optimization in Constrained Robotic Systems}, 
  year={2023},
  volume={39},
  number={1},
  pages={183-202},
  ISSN={1941-0468},
  month={Feb},}

@Book{boumal2023intromanifolds,
  title     = {An introduction to optimization on smooth manifolds},
  author    = {Boumal, Nicolas},
  publisher = {Cambridge University Press},
  year      = {2023},
}

@article{Michael_Ross_review,
title = {A review of pseudospectral optimal control: From theory to flight},
journal = {Annual Reviews in Control},
volume = {36},
number = {2},
pages = {182-197},
year = {2012},
issn = {1367-5788},
author = {I. Michael Ross and Mark Karpenko}
}

@article{garg_RadauPS,
author = {Garg, Divya and Patterson, Michael and Francolin, Camila and Darby, Christopher and Huntington, Geoffrey and Hager, William and Rao, Anil},
year = {2011},
month = {06},
pages = {335-358},
title = {Direct Trajectory Optimization and Costate Estimation of General Optimal Control Problems Using a Radau Pseudospectral Method},
volume = {49},
journal = {Computational Optimization and Applications},
}

@INPROCEEDINGS{mao_scvx,
  author={Mao, Yuanqi and Szmuk, Michael and A\c{c}\i{}kme\c{s}e, Beh\c{c}et},
  booktitle={2016 IEEE 55th Conference on Decision and Control (CDC)}, 
  title={Successive Convexification of Non-Convex Optimal Control Problems and Its Convergence Properties}, 
  year={2016},
  volume={},
  number={},
  pages={3636-3641},
}

@article{szmuk_stc,
author = {Szmuk, Michael and Reynolds, Taylor P. and A\c{c}\i{}kme\c{s}e, Beh\c{c}et},
title = {Successive Convexification for Real-Time Six-Degree-of-Freedom Powered Descent Guidance with State-Triggered Constraints},
journal = {Journal of Guidance, Control, and Dynamics},
volume = {43},
number = {8},
pages = {1399-1413},
year = {2020},
}
                                                   
\appendix
\section{Quaternion retraction: closed-form terms and transports}\label{app:quat_terms}
In this appendix, the closed-form expressions of the retraction-induced 
terms used in Section~6 are summarized for the unit-quaternion manifold $\mathcal Q$.
All expressions below are given in the minimal coordinates
$v\in\mathbb R^3$.
Let $J_r(\phi)\in\mathbb R^{3\times 3}$ denote the right Jacobian. 

\subsection{Closed-form of $C$, $S$, and $\mathcal E$}
For the quaternion kinematics and a reference
$(\bar q,\bar\omega_B)$, the retraction-induced terms in \eqref{eq:ABCS_defs} satisfy
\begin{align}
  C_i^q &= 0,~~S_{i,\rho_i^q}^q(v) = [\rho_i^q]_\times v,~~\mathcal E_i^q(v) = 0,
\end{align}
where $\rho_i^q\in\mathbb R^3$ is the defect.

\subsection{Closed-form of the transport matrices $[{T}]_{ik}^{h}$}\label{app:Tmatrix}
In the numerical example, $\mathcal M$ comprises Euclidean components and $\mathcal Q$; thus, $[{T}]_{ik}^{h}$ is block-diagonal with an identity block and a quaternion block $J_r(\phi)$, satisfying
\begin{align}
  [{T}]_{ik}^{h}& = \mathrm{diag}\!\left(I,\ J_r(\phi_{ki}^{h})\right),\\
  \begin{split}
      J_r(\phi)&=
  I_3
  - \frac{1-\cos\theta}{\theta^2}\,[\phi]_\times
  + \frac{\theta-\sin\theta}{\theta^3}\,[\phi]_\times^2.
  \end{split}\\
 \phi_{ki}^{h}&\coloneqq 2\,\mathrm{Log}\!\left((\bar q_k^{h})^{-1}\otimes \bar q_i^{h}\right)\in\mathbb R^3, \mathrm{where\;} \theta =\|\phi\|.
\end{align}
                                                          
\end{document}